%% file: chaos_2012.tex
\documentclass[11pt]{amsart}
\usepackage{epsfig,color}

\newtheorem{defin}{Definition}

\newcommand{\A}{{\mathcal A}}

\newcommand{\E}{{\mathcal E}}

\newcommand{\HH}{{\mathbb H}}

\newcommand{\mm}{{\mathcal M}}

\newcommand{\N}{{\mathbb N}}\newcommand{\NN}{{\mathbb N}^2}

\newcommand{\pp}{{\mathcal P}}

\newcommand{\R}{{\mathbb R}}\newcommand{\RR}{{\mathbb R}^2}

\newcommand{\SSS}{{\mathbb S}}

\newcommand{\ttt}{{\mathcal T}}

\newcommand{\x}{{\mathcal X}}

\newcommand{\ZZ}{{\mathbb Z}^2}

\newcommand{\bo}{\partial} 

\newcommand{\const}{\mbox{const}} 

\newcommand{\mac}{\text{SL}(2,\R)}   

\newcommand{\tP}{\tilde{P}}

\newcommand{\al}{\alpha}

\newcommand{\ga}{\gamma}\newcommand{\Ga}{\Gamma}

\newcommand{\ep}{\varepsilon}

\newcommand{\Om}{\Omega}

\newcommand{\vp}{\varphi}
\newcommand{\het}{\theta}

\begin{document}

\bibliographystyle{plain}

\title[Problems on Billiards]
{Billiard Dynamics: An Updated Survey with the Emphasis on Open
Problems}
\author{Eugene Gutkin}
\address{Nicolaus Copernicus University, Chopina 12/18, Torun 87-100;
IM PAN, Sniadeckich 8, Warszawa 10, Poland}
\email{gutkin@mat.umk.pl,gutkin@impan.pl}
\date{\today}

\begin{abstract}
This is an updated and expanded version of our earlier survey
article \cite{Gut5}. Section $\S 1$ introduces the subject matter.
Sections $\S 2 - \S 4$ expose the basic material following the
paradigm of elliptic, hyperbolic and parabolic billiard dynamics.
In section $\S 5$ we report on the recent work pertaining to the
problems and conjectures exposed in the survey \cite{Gut5}.
Besides, in section $\S 5$ we formulate a few additional problems
and conjectures. The bibliography has been updated and
considerably expanded.
\end{abstract}

\maketitle

\tableofcontents

\section{Introduction} \label{intro}
Billiard dynamics broadly understood is the geodesic flow on a
Riemannian manifold with a boundary. But even this very general
framework is not broad enough, e.g., for applications in physics.
In these applications the manifold in question is the
configuration space of a physical system. Often, it is a {\em
manifold with corners} and {\em singularities}. Some physics
models lead to the Finsler billiard \cite{GT}: The  manifold in
question is not Riemannian; it is Finslerian. The simplest
examples of Riemannian manifolds with corners are plane polygons,
and some basic physical models yield the billiard on triangles
\cite{Kolan,GuTr,Glashow,Gut3}

The configuration space of the famous gas of elastic balls
\cite{Sinai-,Sas1} is structured combinatorially like a euclidean
polyhedron of a huge number of dimensions. In fact, this
configuration space is much more complicated,  because the
polyhedron is not flat. The mathematical investigation of this
system produced the celebrated Boltzmann Ergodic Hypothesis. After
Sinai's seminal papers \cite{Sinai-,Sinai}, a modified version of
the original conjecture became known as the Boltzmann-Sinai
Hypothesis.

However, the bulk of our exposition is restricted to the billiard
in a bounded planar domain with piecewise smooth boundary. The
reason is threefold. First of all, this setting allows us to avoid
lengthy preliminaries and cumbersome formalism: It immediately
leads to qualitative mathematical questions. (This was also the
opinion of G. D. Birkhoff \cite{Birk}.) Second, there are basic
physical models that correspond to planar billiards \cite{Gut3}.
Third, and most important, there are fundamental problems on the
plane billiard that are still open. The problems are indeed
fundamental: They concern the main features of these dynamical
systems.

In the body of the paper  we introduce several open problems of
billiard dynamics. Our choice of the questions is motivated partly
by the personal taste and partly by the simplicity of formulation.
We review the preliminaries, discuss the motivation, and outline
possible angles of attack. We also point out partial results and
other evidence toward the answer. Formally, the exposition is
self-contained, but the reader may want to consult the literature
\cite{KH,HK,KozTres,Tabach,Chernov4,Gut2,Gut3}.

For obvious reasons, we will call the planar domain in question
the {\em billiard table}. Its geometric shape determines the
qualitative character of the motion. Historically, three classes
of shapes have mostly attracted attention. First, it is the class
of smooth and strictly convex billiard tables. For several
reasons, the corresponding billiard dynamics is called elliptic.
Second, it is the piecewise concave and piecewise smooth billiard
tables. The corresponding dynamics is hyperbolic.\footnote{There
are also convex billiard tables that yield hyperbolic dynamics.
See section~~\ref{hyperbolic}.} Billiard tables of the third class
are the polygons. The corresponding dynamics is parabolic. The
three types of the billiard are exposed in \S~\ref{elliptic} --
\S~\ref{hyperbolic} -- \S~\ref{parabolic} repectively.

%

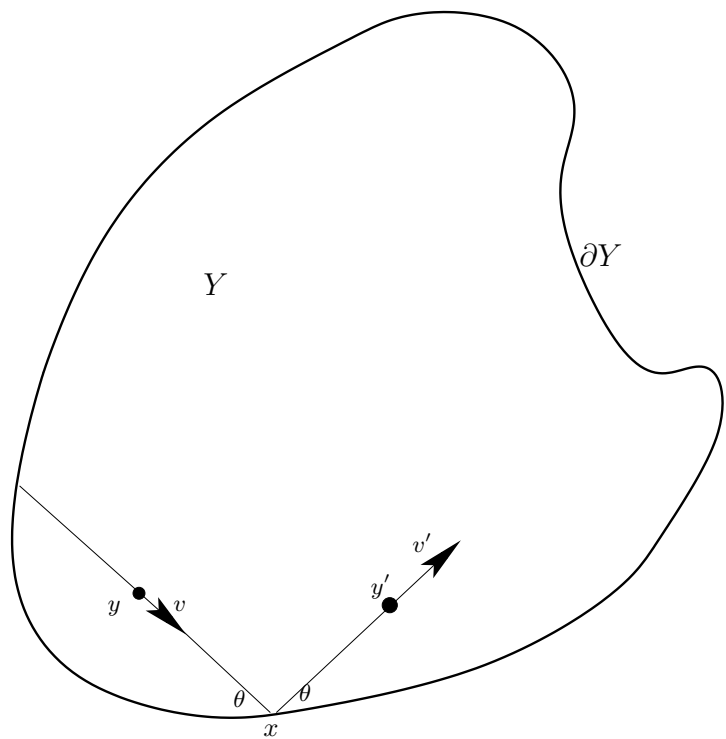
\begin{figure}[htbp]
\begin{center}
\input{billi1.pstex_t}
\caption{The billiard flow.}
    \label{fig1}
\end{center}
\end{figure}

In the rest of the introduction we describe the basic notation and
the terminology. Let $Y\subset\RR$ be a compact, connected
billiard table. See Figure~\ref{fig1}. Its boundary $\bo Y$ is a
finite union of $C^1$ curves. It may have several connected
components. The {\em billiard flow} on $Y$ is modelled on the
motion of a material point: The ``particle" or the ``billiard
ball". At each time instant, the state of the system is determined
by the position of the ball, $y\in Y$, and its velocity, a unit
vector $v\in\RR$. (It suffices to consider the motion with the
unit speed.) The ball rolls along the ray emanating from $y$, in
the direction $v$. At the instant the ball reaches $\partial Y$,
its direction changes. Let $x\in\partial Y$ be the point in
question, and let $v'$ be the new direction. The transformation,
$v\mapsto v'$, is the {\em orthogonal reflection} about the
tangent line to $\partial Y$ at $x$. The vector $v'$ is directed
inward, and the ball keeps rolling.

These rules define: 1) The {\em phase space} $\Psi$ of the
billiard flow, as the quotient of $Y\times S^1$ by the
identification $(x,v)=(x,v')$ above; 2) The billiard flow
$T^t:\Psi\to\Psi$. If $Y$ is simply connected, and $\partial Y$ is
$C^1$, then $\Psi$ is homeomorphic to the three-dimensional
sphere.\footnote{We are not aware of any uses of this observation
in the billiard literature.} In any way, $\dim\Psi=3$, and the
reader may think of $\Psi$ as the set of pairs $(y,v)$, such that
$v$ is directed inward.

A few remarks are in order. The rules defining the billiard flow
stem from the assumptions that the billiard motion is
frictionless, and that the boundary of the billiard table is
perfectly elastic. The orthogonal reflection rule $v\mapsto v'$
insures that billiard orbits are the local minimizers of the
distance functional. (This property extends to  the Finsler
billiard \cite{GT}.) The reflection rule is not defined at the
corners of the boundary. The standard convention is to ``stop the
ball" when it reaches a corner. Thus, if $\bo Y$ is not $C^1$,
then there are billiard orbits that are not defined for all times.
Their union has zero volume with respect to the {\em Liouville
measure} defined below.

Set $X=\bo Y$, and endow it with the positive orientation.
Choosing a reference point on each connected component, and using
the arc length parameter, we identify $X$ with the disjoint union
of $k\ge 1$ circles. In this paper, with the exception of
\S~\ref{parabolic}, $k=1$. The set $\Phi\subset\Psi$ given by the
condition $y\in\partial Y$ is a cross-section for the billiard
flow. The Poincar\'e mapping $\vp:\Phi\to\Phi$ is the {\em
billiard map} and $\Phi$ is its phase space. The terminology is
due to G. D. Birkhoff who championed the ``billiard ball problem"
\cite{Birk}. Let $x$ be the arc length parameter on $X$. For
$(x,v)\in\Phi$ let $\theta$ be the angle between $v$ and the
positive tangent to $\partial Y$ at $x$. Then $0\le\theta\le\pi$,
where $0$ and $\pi$ correspond to the forward and the backward
tangential directions respectively. This coordinate system fails
at the corners of $\bo Y$. If $\bo Y$ is $C^1$, then
$\Phi=X\times[0,\pi]$. We will use the notation
$\vp(x,\theta)=(x_1,\theta_1)$.

%

\begin{figure}[htbp]
\begin{center}
\input{billi2.pstex_t}
\caption{The billiard map.}
    \label{fig2}
\end{center}
\end{figure}
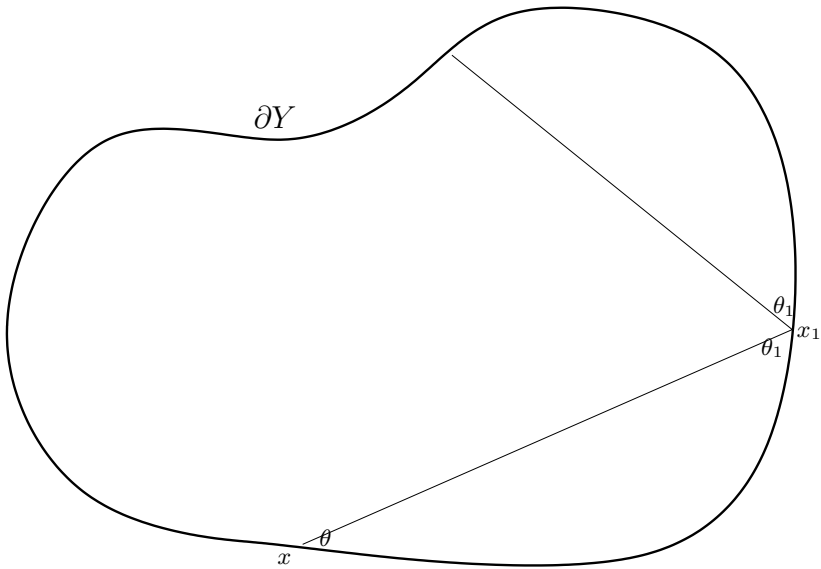

Let $p,q$ be the euclidean coordinates in $\RR$, and let $0\le
\alpha < 2\pi$ be the angle coordinate on the unit circle. The
Liouville measure on $\Psi$ has the density $d\nu=dpdqd\alpha$. It
is invariant under the billiard flow. The {\em induced Liouville
measure} $\mu$ on $\Phi$ is invariant under the billiard map, and
has the density $d\mu=\sin\theta dxd\theta$. Both measures are
finite. Straightforward computations yield
\begin{equation} \label{volumes-eq}                 
\nu(\Psi)=2\pi\mbox{Area}(Y),\ \mu(\Phi)=2\,\mbox{Length}(\partial Y).
\end{equation}
\section{Smooth, strictly convex billiard: elliptic dynamics} \label{elliptic}
The first deep investigation of this framework is due to G. D.
Birkhoff \cite{Birk}. For this reason, it is often called the {\em
Birkhoff billiard}. The billiard map is an {\em area preserving
twist map} \cite{KH}.  An {\em invariant circle} is a
$\vp$-invariant curve $\Gamma\subset\Phi$ which is a
noncontractible topological circle. Recall that $\Phi$ is a
topological annulus. Both components of $\partial\Phi$ are the
trivial invariant circles. From the geometric optics viewpoint,
$\Phi$ is the space of light rays (i.e., directed lines), and $Y$
is a room whose walls are the perfect mirrors. Then
$\Gamma\subset\Phi$ is a one-parameter family of light rays in
$Y$, and its {\em envelope} $F(\Gamma)$ is the set of focusing
points of light rays in this family. Note that $F(\Gamma)$ is not
a subset of $Y$, in general. For instance, if $Y$ is an ellipse,
then there are invariant curves $\Gamma$ such that $F(\Gamma)$ are
confocal hyperbolas.

\begin{figure}[htbp]
\begin{center}
\input{billi3.pstex_t}
\caption{The phase space of the billiard map in an ellipse.}
    \label{fig3}
\end{center}
\end{figure}
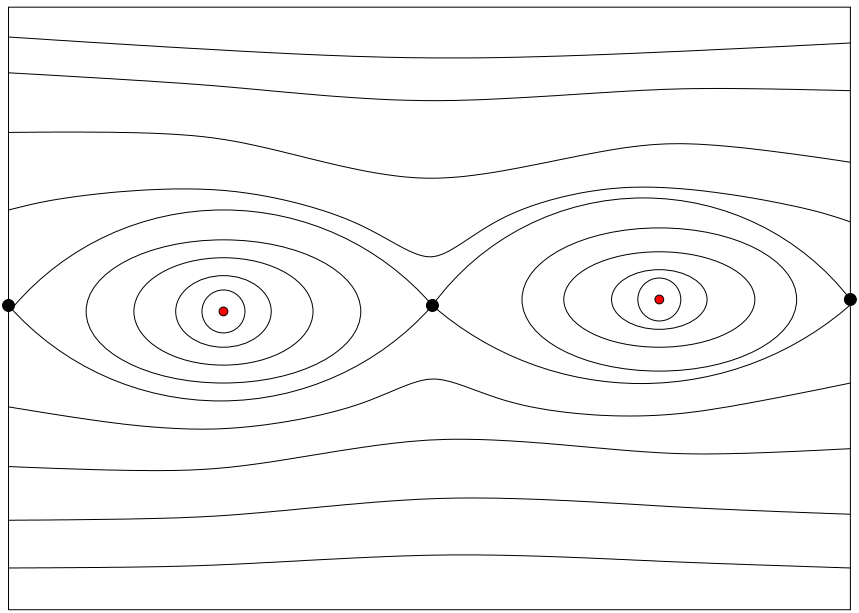

Let $\Gamma$ be an invariant circle, and let $\gamma=F(\Gamma)$.
Then $\gamma\subset\mbox{Int}(Y)$ \cite{GutKat}. These curves are
the {\em caustics} of the billiard table. When $\bo Y$ is an
ellipse, the caustics are the confocal ellipses. Their union is
the region $Y\setminus[ff']$, where $f,f'$ are the foci of $\in
Y$. If $Y$ is not a disc, the invariant circles fill out a region,
$C(\Phi)\subset\Phi$, whose complement looks like a pair of
``eyes". See Figure~\ref{fig3}.


\begin{defin}  \label{integrable-def}  
A billiard table $Y$ is {\em integrable} if the set of invariant
circles has nonempty interior.
\end{defin}
The most famous open question  about caustics is known as the {\em
Birkhoff conjecture}. It first appeared in print in a paper by
Poritsky \cite{Porit}, several years after Birkhoff's
death.\footnote{In the introduction to \cite{Porit}, the author
says that many years ago, when he was a doctoral student of
Birkhoff, his advisor communicated the conjecture to him.}

\medskip

\noindent {\bf Problem 1 (Birkhoff conjecture)}. Ellipses are the only integrable billiard tables.\\

A disc is a degenerate ellipse, with $f=f'$. The preceding
analysis applies, and the invariant circles fill out all of the
phase space. M. Bialy proved the converse: If all of $\Phi(Y)$ is
foliated by invariant circles, then $Y$ is a disc \cite{Bialy}.
See  \cite{Bialy1} for an extension of this theorem to the
surfaces of arbitrary constant curvature. We refer the reader to
section \S~~\ref{update} for elaborations and updates on the
Birkhoff conjecture.

\medskip

Let $Y$ be any oval. If $X=\partial Y$ is sufficiently smooth, and
its curvature is strictly positive, then the invariant circles
fill out a set of positive measure. This was proved by V. Lazutkin
under the assumption that $\partial Y$ was of class $C^{333}$
\cite{Laz}. Lazutkin's proof crucially uses a famous theorem of J.
Moser \cite{Moser1}.\footnote{This is a seminal paper on the {\em
KAM theory}. See also \cite{Moser2}. The name KAM stands for
Kolmogorov, Arnold and Moser. See \cite{deLla} for an exposition.}
The number $333$ is chosen in order to satisfy the assumptions in
\cite{Moser1}. The required smoothness was eventually lowered to
$C^6$~~\cite{Douad}. By a theorem of J. Mather \cite{Mather1}, the
positive curvature condition is necessary for the existence of
caustics.

An invariant noncontractible topological annulus,
$\Omega\subset\Phi$, whose interior contains no invariant circles,
is a {\em Birkhoff instability region}. This is a special case of
an important concept for area preserving twist maps \cite{KH}.
Assume the Birkhoff conjecture, and let $Y$ be a non-elliptical
billiard table. Then $\Phi$ contains Birkhoff instability regions.
The dynamics in an instability region has positive {\em
topological entropy} \cite{Angen}. Hence, the Birkhoff conjecture
implies that any non-elliptical billiard has positive topological
entropy. By the {\em (metric) entropy} of a billiard we will mean
the entropy of the Liouville measure. The only examples of convex
billiard tables with positive entropy are the Bunimovich stadium
\cite{Buni+} and its generalizations. These billiard tables are
not strictly convex, and their boundary is only $C^1$. The
corresponding billiard dynamics is
hyperbolic. See \S~\ref{hyperbolic}. This leads to our next open question.\\

\noindent {\bf Problem 2}. a) Construct a
strictly convex $C^1$-smooth billiard table with positive entropy.
b) Construct a convex $C^2$-smooth billiard table with positive entropy.\\

Using an ingenuous variational argument, Birkhoff proved the
existence of certain periodic billiard orbits \cite{Birk}. His
approach extends to area preserving twist maps, and thus yields a
more general result on periodic orbits of these dynamical systems
\cite{KH}. In the billiard framework the relevant considerations
are especially transparent. A periodic orbit of period $q$
corresponds to an (oriented) closed polygon with $q$ sides,
inscribed in $Y$, and satisfying the obvious condition on the
angles it makes with $\partial Y$. Birkhoff called these the {\em
harmonic polygons}. Vice versa, any oriented harmonic $q$-gon $P$
determines a periodic orbit of period $q$. Let $1\le p <q$ be the
number of times the pencil tracing $P$ goes around $\partial Y$.
The ratio $0 < p/q < 1$ is the {\em rotation number} of a periodic
orbit. Fix a pair $1\le p <q$, with $p$ and $q$ relatively prime.
Let $X(p,q)$ be the set of all inscribed $q$-gons that go $p$
times around $\partial Y$. The space $X(p,q)$ is a manifold with
corners. For $P\in X(p,q)$ let $f(P)$ be the circumference of $P$.
Then harmonic polygons are the critical points of the function
$f:X(p,q)\to\R$. Birkhoff proved that $f$ has at least two
distinct critical points. One of them delivers the maximum, and
the other a minimax to the circumference. The corresponding
periodic billiard orbits are the {\em Birkhoff periodic orbits}
with the rotation number $p/q$.

By way of example, we take the rotation number $1/2$. Then the
maximal Birkhoff orbit yields the {\em diameter of $Y$}. The
minimax orbit corresponds to the {\em width of $Y$}. When the
diameter and the width of $Y$ are equal, the boundary $\partial Y$
is a {\em curve of constant width}; then we have a one-parameter
family of periodic orbits with the rotation number $1/2$. They
fill out the ``equator" of $\Phi$. There are other examples of
ovals with one-parameter families of periodic orbits having the
same length and the same rotation number. See \cite{Innami} and
\cite{Gut+,Gut10} for different approaches.


One of the basic characteristics of a dynamical system is the {\em
growth rate of the number of periodic points}. In order to talk
about it, we need a {\em counting function}. The standard counting
function $f_Y(n)$ for the billiard map is the number of periodic
points of the period at most $n$. (See \S~\ref{hyperbolic} and
\S~\ref{parabolic}  for other examples.) The set of periodic
points is partitioned into periodic orbits, and let $F_Y(n)$ be
the number of periodic orbits of period at most $n$. Birkhoff's
theorem bounds $F_Y(n)$ from below by the number of relatively
prime pairs $1\le p < q \le n$. This implies a universal cubic
lower bound $f_Y(n)\ge cn^3$.  See, e. g., \cite{Hardy}.

\medskip

Since an oval may have infinitely many periodic points of the same
period, there is no universal upper bound on $f_Y(n)$. The size of
a measurable set is naturally estimated by its measure. Let
${\mathcal P}\subset\Phi$ (resp. ${\mathcal P}_n\subset\Phi$) be
the set of periodic points (resp. periodic points of period $n$).
For example, if $Y$ is a table of constant width, then ${\mathcal
P}_2\subset\Phi$ is the equator. Although it is infinite,
$\mu({\mathcal P}_2)=0$. Since ${\mathcal
P}=\cup_{n=2}^{\infty}{\mathcal P}_n$, a disjoint union,
$\mu({\mathcal P})=\sum_{n=2}^{\infty}\mu({\mathcal P}_n)$. Thus,
$\mu({\mathcal P})=0$ iff $\mu({\mathcal P}_n)=0$ for all
$n=2,3,4,\dots$.

The famous {\em Weyl formula} gives the leading term and the error
estimate for the spectral asymptotics  of the Laplace operator
(with either Dirichlet or Neumann boundary conditions) in a
bounded domain of the euclidean space (of any number of
dimensions). The (also famous) {\em Weyl conjecture} predicts the
second term of the asymptotic series \cite{Weyl}. A theorem of V.
Ivrii \cite{Ivrii} establishes the Weyl conjecture for a euclidean
domain under the assumption that the set of periodic billiard
orbits has measure zero.\footnote{A more general formula for the
spectral asymptotics of the Laplacean, due to Safarov and
Vassiliev, contains a term accounting for periodic orbits
\cite{SaVa}. If periodic points yield a set of measure zero, this
term vanishes.}

Ivrii conjectured that the assumption $\mu({\mathcal P})=0$ was
superfluous: It should hold for any euclidean domain with a smooth
boundary. Members of the Sinai's dynamics seminar in Moscow
promised to him in 1980 to prove the desideratum in a few days ...
The question is still open. Problem 3 below states the conjecture
for plane domains.\\

\noindent{\bf Problem 3 (Ivrii conjecture)}. Let $Y$ be a
piecewise smooth billiard table.
i) Prove that $\mu({\mathcal P})=0$. ii) Prove that $\mu({\mathcal P}_n)=0$ for all $n$.\\

Although Problem 3 concerns arbitrary billiard tables, it is
especially challenging for the Birkhoff billiard, hence we have
put the problem into this section. It is convenient to designate
by, say, $I_n$ the claim $\mu({\mathcal P}_n)=0$. Thus, Ivrii
conjecture amounts to proving $I_n$ for all $n\ge 2$. Claim $I_2$
is obvious. Proposition $I_3$ is a theorem of M. Rychlik
\cite{Rych}. His proof depends on a formal identity, verified
using Maple. L. Stojanov simplified the proof, and eliminated the
computer verification \cite{Stoj2}. Ya. Vorobets gave an
independent proof \cite{Vor1}. His argument applies to higher
dimensional billiards as well. M. Wojtkowski \cite{Wojt3} obtained
Rychlik's theorem as an application of the {\em mirror equation}
of the geometric optics and the {\em isoperimetric inequality}.
See \cite{GutKat} for other applications.

Ivrii's conjecture is known to hold in many special cases, e.g.,
for hyperbolic and parabolic billiard tables. See
\S~\ref{hyperbolic} and \S~\ref{parabolic}. It holds for billiard
tables with real analytic boundary \cite{SaVa}. For the generic
billiard table the sets ${\mathcal P}_n$ are finite for all $n$
\cite{PetStoj}. The billiard map for a Birkhoff billiard table is
an area preserving twist map. However, there are smooth area
preserving twist maps such that $\mu({\mathcal P})>0$. Thus,
Ivrii's conjecture is really about the billiard map!

Recently Glutsyuk and Kudryashov announced a proof of Proposition
$I_4$ \cite{GluKudr}. See Section~~\ref{elliptic_sub} for further
comments.

\section{Hyperbolic billiard dynamics} \label{hyperbolic}
It is customary to say that a billiard table is hyperbolic if the
associated dynamics is hyperbolic. The dynamics in question may be
the billiard flow or the billiard map or the induced map on a
subset of the phase space. For concreteness, we will call a
billiard table hyperbolic if the corresponding billiard map is
hyperbolic. The modern approach to hyperbolic dynamics crucially
uses the Oseledets multiplicative ergodic theorem \cite{Osel}. See
\cite{KH,HK} for a general introduction into the hyperbolic
dynamics and \cite{KS,Buni,Tabach} (resp. \cite{Chernov4}) for
introductory (resp. thorough) expositions of the hyperbolic
billiard.

The first hyperbolic billiard tables were made from concave arcs.
As a motivation, let us consider the following construction. Let
$P$ be a convex polygon. Replace some of the sides of $P$ by
circular arcs whose centers are sufficiently far from $P$. The
result is a ``curvilinear polygon", $Y$, approximating $P$.
Choosing appropriate center points, we insure that the ``curved
sides" of $Y$ are convex inward. It is not important that they be
circular, as long as they are smooth and convex inward.

\begin{figure}[htbp]
\begin{center}
\input{billi4.pstex_t}
\caption{The Sinai billiard table.}
    \label{fig4}
\end{center}
\end{figure}
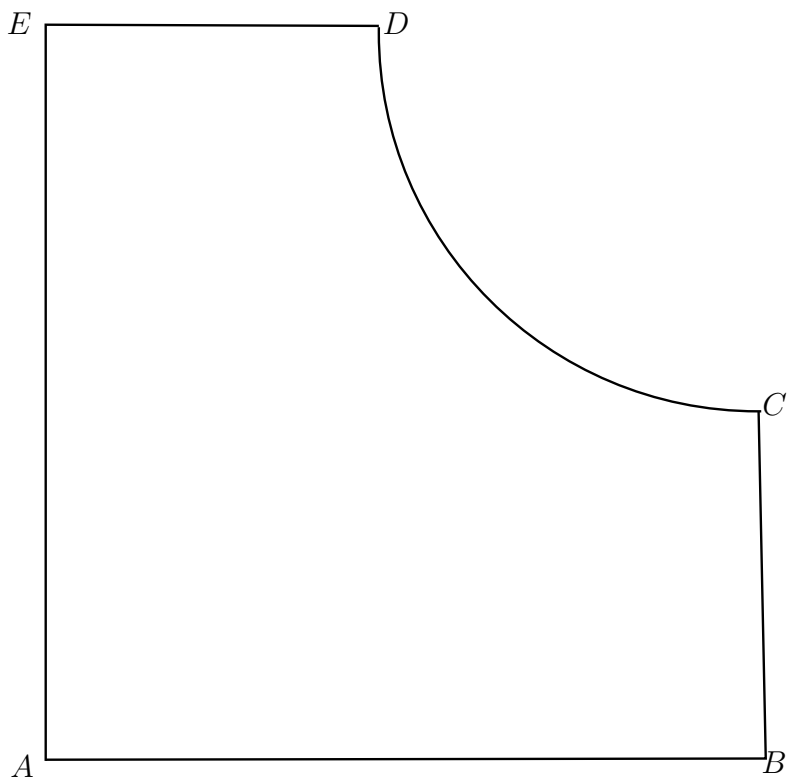

This class of billiard tables arose in the work of Ya. Sinai on
the {\em Boltzmann-Sinai gas} \cite{Sinai-}.\footnote{See the
Appendix by D. Szasz in \cite{Sas1}.} In the {\em Boltzmann gas}
the identical round molecules are confined by a box. Sinai has
replaced the box by periodic boundary conditions. Thus, the
molecules of the Boltzmann-Sinai gas move on a flat torus. In the
``real world", the confining box is three-dimensional and the
number of moving molecules is enormous. In the Sinai
``mathematical caricature", there are only two molecules on a
two-torus. The system reduces to the geodesic flow on a flat torus
with a round hole. Represent the flat torus by the $2\times 2$
square, so that the hole is the central disc of radius $1/2$. By
the four-fold symmetry, the problem reduces to the billiard on the
unit square with the deleted quarter-disc of radius $1/2$,
centered at a vertex. See Figure~\ref{fig4}.


This domain is known as {\em the Sinai
billiard}.\footnote{Unfortunately, there is a fair amount of
confusing terminology in the literature. Mathematicians often use
the expressions like ``Sinai's billiard" or ``the Sinai table" or
``a dispersive billiard" interchangeably. Physicists tend to mean
by ``the Sinai billiard" a special billiard table, although not
necessarily that of Figure 4.} Let now $P$ be a (not necessarily
convex) $n$-gon. Let $Y$ be the region obtained by replacing $1\le
m < n$ (resp. all) of the sides of $P$ by circular arcs,
satisfying the conditions above. Then $Y$ is a {\em
semi-dispersive (resp. dispersive) billiard table}. The circular
arcs (resp. the segments) of $\bo Y$ are its {\em dispersive
(resp. neutral) components}. This terminology extends to the
billiard tables whose dispersive boundary components are smooth,
convex inward curves. These are the (semi)dispersive billiard
tables. The table in Figure 4 has one dispersive and four neutral
boundary components.

In \cite{Sinai} Sinai proved the hyperbolicity of dispersive
billiard tables. After the discovery by L. Bunimovich that the
{\em stadium} and similar billiard tables are hyperbolic
\cite{Buni+}, mathematicians started searching for geometric
criteria of hyperbolicity. The notion of an {\em invariant cone
field} \cite{Wojt+,KB} proved to be very useful.

Denote by $V_z$ the tangent plane to the phase space at
$z\in\Phi$. The differential $\vp_*$ is a linear map from $V_z$ to
$V_{\vp(z)}$. By our convention, a subset of a vector space is a
{\em cone} if it is invariant under multiplications by all
scalars.

\begin{defin} \label{cone-def}                   
A family ${\mathcal C}=\{C_z\subset V_z:z\in\Phi\}$ is an {\em
invariant cone field} if the following conditions are satisfied.
\begin{itemize}
\item 1. The closed cone $C_z$ is defined for almost all $z\in\Phi$,
and the map $z\mapsto C_z$ is measurable.
\item 2. The cone $C_z$ is has nonempty interior.\\
\item 3. We have $\vp_*(C_z)\subset C_{\vp(z)}$.\\
\item 4. There exists $n=n(z)$ such that
$\vp_*^n(C_z)\subset \mbox{int}(C_{\vp^n(z)})$.\\
\end{itemize}
\end{defin}
The hyperbolicity is equivalent to the existence of an invariant
cone field \cite{Wojt+}. Wojtkowski constructed invariant cone
fields for several classes of billiard tables \cite{Wojt1}. In
addition to the dispersive tables and the generalized stadia, he
found invariant cone fields for a wide class of locally strictly
convex tables. Wojtkowski's approach was further extended by
Bunimovich, V. Donnay, and R. Markarian \cite{Chernov4}. Using
these ideas, B. Gutkin, U. Smilanski and the author constructed hyperbolic
billiard tables on surfaces of arbitrary constant curvature~~\cite{GuSm}.\\

\noindent {\bf Problem 4}. Is every semi-dispersive billiard table
hyperbolic?\\

Let $Y$ be a semi-dispersive $n$-gon with only one neutral
component. Let $Y'$ be the reflection of $Y$ about this side, and
set $Z=Y\cup Y'$. Since $Z$ is a dispersive billiard table, it is
hyperbolic. By the reflection symmetry, the table $Y$ is also
hyperbolic. In special cases, the reflection trick yields the
hyperbolicity of semi-dispersive $n$-gons with $m<n-1$ dispersive
components. For instance, let $P$ be a triangle with an angle
$\pi/n$. Let $Y$ be the semi-dispersive triangle, whose only
dispersive component is located opposite the $\pi/n$ angle.
Reflecting $Y$ successively $2n$ times, we obtain a dispersive
billiard table, $Z$. Thus, $Z$ is hyperbolic. By the symmetry, the
table $Y$ is hyperbolic as well. A suitable generalization of the
reflection trick will work if $P$ is a {\em rational polygon}. See
\S~\ref{parabolic}. The special case $m=1$ of Problem 4 is closely
related to Problem 9 of \S~\ref{parabolic}.

Dispersive billiard tables are ergodic \cite{BSC2}. There are
examples of hyperbolic, but nonergodic billiard tables
\cite{Wojt1}. The consensus is that a typical hyperbolic billiard
is ergodic. For instance, the stadium and its relatives are
ergodic \cite{Szasz}. There are no examples of strictly convex
hyperbolic billiard tables. See Problem 2.

For the rest of this section, we consider only dispersive billiard
tables. Referring the reader to
\cite{BSC1,BSC2,Chernov1,Sas1,Chernov2,LaiSun} for a discussion of
their chaotic properties and to open questions about, e.g., the
{\em decay of correlations}, we concentrate on the statistics of
periodic orbits in hyperbolic billiards. The set of periodic
points of any period is finite; let $f_Y(n)$ be the number of
periodic points, whose period is less than or equal to $n$. The
asymptotics of $f_Y(n)$, as $n\to\infty$, is an important
dynamical characteristic. By theorems of Stojanov and Chernov
\cite{Stoj,BSC1}, there are $0<h_-<h_+<\infty$ such that
\begin{equation} \label{lower-upper-eq} 
0<h_-\le\liminf_{n\rightarrow\infty}\frac{\log f_Y(n)}{n}\le
\limsup_{n\rightarrow\infty}\frac{\log f_Y(n)}{n}\le h_+<\infty.
\end{equation}
The following two problems were contributed by N. Chernov.\\

\noindent {\bf Problem 5}. Does the limit
\begin{equation} \label{asympt-eq}  
h=\lim_{n\rightarrow\infty}\frac{\log f_Y(n)}{n}
\end{equation}
exist?\\

\noindent {\bf Problem 6}. If the limit in
equation~(\ref{asympt-eq}) exists, is $0<h<\infty$
the {\em topological entropy} of the billiard map?\\

Problems 5 and 6 fit into the general relationship between the
distribution of periodic points and the topological entropy
\cite{Kat2}. However, the singularities, which constitute the
paramount feature of billiard dynamics, preclude the applicability
of smooth ergodic theory. Other techniques have to be developed
\cite{Chernov3,GuHa95,GuHa}.

\section{Polygonal billiard: parabolic dynamics} \label{parabolic}
The polygon $P$ that serves as a billiard table is not required to
be convex or simply connected. It may also have  {\em barriers},
i. e., obstacles without interior. It is {\em rational} if the
angles between its sides are of the form $\pi m/n$. Let $N=N(P)$
be the least common denominator of these rational numbers. A
classical construction associates with $P$ a closed surface
$S=S(P)$ tiled by $2N$ copies of $P$. The surface $S$ has a finite
number of cone points; the cone angles are integer multiples of
$2\pi$. Suppose that $P$ is a {\em simple polygon},\footnote{A
polygon $P$ is simple if $\bo P$ is connected.} and let
$m_i\pi/n_i,1\le i \le p,$ be its angles. The genus of $S(P)$
satisfies \cite{Gut2}
\begin{equation} \label{genus-eq}  
g(S(P)) = 1+\frac{N}{2}\sum_{i=1}^p\frac{m_i-1}{n_i}.
\end{equation}
Equation~(\ref{genus-eq}) implies that $S(P)$ is a torus if and
only if $P$ tiles the plane under reflections. The billiard in $P$
is essentially equivalent to the geodesic flow on $S(P)$. This
observation was first exploited by A. Katok and A. Zemlyakov
\cite{KZ}, and $S(P)$ is often called the ``Katok-Zemlyakov
surface". However, the construction has been in the literature (at
least) since the early 20-th century \cite{Stackel,Fox}. We refer
to the surveys \cite{Gut2,Gut3,Smillie,MasTab,Tabach} for
extensive background material.

Surfaces $S(P)$ are examples of {\em translation surfaces}, which
are of independent interest \cite{GuJ}. From the viewpoint of
classical analysis, a translation surface is a closed Riemann
surface with a holomorphic linear differential. Using holomorphic
quadratic (as opposed to linear) differentials, we arrive at the
notion of {\em half-translation surfaces} \cite{Guj,GuJ}. Billiard
orbits on a polygon become geodesics on the corresponding
translation (or the half-translation) surface. Since billiard
orbits change directions at every reflection, the notion of the
direction of an orbit is not well defined. Geodesics on a
translation surface, on the contrary, do not change their
directions. This yields a technical advantage of translation
surfaces over polygons \cite{KZ}. The crucial advantage comes,
however, from the natural action of the group $\mac$ on
translation surfaces \cite{KMS,Mas1,Mas2,Veech2,Smillie}. See
section~~\ref{update} for elaborations.

The geodesic flow of any translation surface, $S$, decomposes into
the one-parameter family of {\em directional flows}
$b_{\theta}^t,\ 0\le \theta < 2\pi$. The flow $b_{\theta}^t$ is
identified with the {\em linear flow on $S$ in direction
$\theta$}. The Lebesgue measure on $S$ is preserved by every
$b_{\theta}^t$. Thus, not only is the billiard flow of a rational
polygon not ergodic, it decomposes as a one-parameter family of
{\em directional billiard flows}. Let $S$ be an arbitrary
translation surface. A theorem of Kerckhoff, Masur, and Smillie
\cite{KMS} says that the flows $b_{\theta}^t$ are uniquely ergodic
for Lebesgue almost all $\theta$. In particular, the directional
billiard flow of a rational polygon is ergodic for almost every
direction. The set ${\mathcal N}(S)\subset[0,2\pi)$ of
non-uniquely ergodic directions has positive Hausdorff dimension
for the typical translation surface \cite{MaSmi}; for particular
classes of rational polygons and translation surfaces the sets
${\mathcal N}(S)$ are countably infinite
\cite{Gut1,Veech2,Yitwah}. We point out that a typical translation
surface does not correspond to any polygon, which illustrates the
limitations of this relationship for the study of polygonal
billiard. See section~~\ref{update} for elaborations on the
polygonal billiard and translation surfaces.\\


Much less is known about the billiard in irrational (i. e.,
arbitrary) polygons. Denote by $\ttt(n)$ the moduli space of
simple euclidean $n$-gons. Since the billiard dynamics is
invariant under scaling, in $\ttt(n)$ we identify polygons that
coincide up to scaling. The space $\ttt(n)$ is a finite union of
components which correspond to particular combinatorial data. We
will refer to them as the {\em combinatorial type components}.
Each component is homeomorphic to a relatively compact set of the
maximal dimension in a Euclidean space. Let $\lambda$ be the
probability measure on $\ttt(n)$, such that its restrictions to
the combinatorial type components are the corresponding Lebesgue
measures. For instance, the space $\ttt(3)\subset\R^2$ is given by
$\ttt(3)=\{(\alpha,\beta): 0< \alpha \le\beta <\pi/2\}$. Thus,
$\ttt(3)$ itself is a plane triangle. By a theorem in \cite{KMS},
the set $\E(n)\subset\ttt(n)$ of ergodic $n$-gons is residual in
the sense of Baire category \cite{Oxtoby}.\\

\noindent {\bf Problem 7}.  Is $\lambda(\E(n))>0$ ?\\

\begin{figure}[htbp]
\begin{center}
\input{billi5.pstex_t}
\caption{Three perfectely elastic particles on the circle.}
    \label{fig5}
\end{center}
\end{figure}
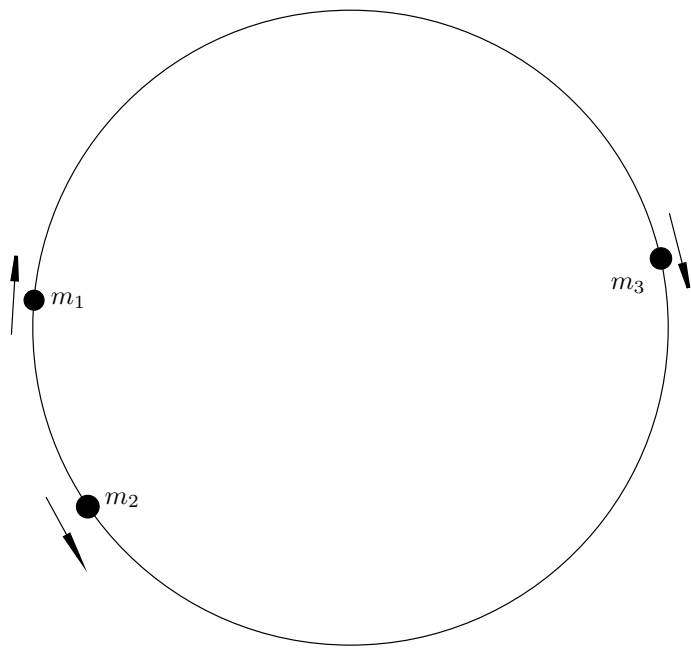

The case of $n=3$ is especially interesting, since the mechanical
system of three elastic point masses moving   on a circle (see
Figure~\ref{fig5}) leads to the billiard in an acute triangle
\cite{Glashow,Casati2}. Let $m_1,m_2,m_3$ be the masses. Then the
angles of the corresponding triangle $\Delta(m_1,m_2,m_3)$ satisfy
\begin{equation}  \label{angles-eq}  
\tan\alpha_i=m_i\sqrt{\frac{m_1+m_2+m_3}{m_1m_2m_3}}.
\end{equation}
We point out that the rationality of the triangle corresponding to
a mechanical system of point masses does not have any obvious
physical meaning. In the limit, when $m_3\rightarrow\infty$, we
obtain the physical system of two elastic particles on an
interval. The limit of $\Delta(m_1,m_2,m_3)$ is the right triangle
whose angles satisfy $\tan\alpha_1=\sqrt{m_1/m_2},
\tan\alpha_2=\sqrt{m_2/m_1}$.

%

Let $P$ be an irrational polygon. Let $\al_1,\dots,\al_k$ be its
angles. If the numbers $\al_i/\pi$ simultaneously admit a certain
super-exponentially fast rational approximation, then $P$ is
ergodic \cite{Vorob}. This remarkable theorem yields explicit
examples of ergodic polygons. However, it does not help with the
above problem. There is some numerical evidence that irrational
polygons are ergodic and have other stochastic properties
\cite{Casati1,Casati2}. So far, there are no theorems confirming
or precluding this.\\

\noindent {\bf Problem 8}. Give an example of an irrational but
nonergodic polygon.\\

Let $P$ be an arbitrary $n$-gon, and let $a_1,\dots,a_n$ be its
sides. For $1\le i \le n$ let $\Phi_i\subset\Phi$ be the set of
elements whose base points belong to the side $a_i$. Then
$\Phi=\cup_{i=1}^n\Phi_i$, a disjoint decomposition. By
equation~(\ref{volumes-eq}), $\mu(\Phi_i)=2\,\mbox{Length}(a_i)$.
The next question/conjecture concerns the structure of invariant
sets in the phase space of a nonergodic polygon. (Compare with
Problem 4 in \S~\ref{hyperbolic}.) By \cite{KMS}, the conjecture holds for rational polygons.\\

\noindent {\bf Problem/Conjecture 9}. Let $P$ be an irrational
$n$-gon, and let $M\subset\Phi$ be an invariant set of positive
measure. If $\Phi_i\subset M$ for some $1\le i \le n$, then $M=\Phi$.\\

The subject of periodic billiard orbits in polygons requires only
elementary euclidean geometry, and has immediate applications to
physics. For instance, let $\Delta(m_1,m_2,m_3)$ be the acute
triangle corresponding to the system of three elastic point masses
equation~(\ref{angles-eq}). Periodic billiard orbits in
$\Delta(m_1,m_2,m_3)$ correspond to the periodic motions of this
mechanical system. \footnote{Another connection between billiards
and physics arises in the study of mechanical linkages
\cite{Sossinsky}.} Ironically, periodic orbits in polygons
turned out to be especially elusive.\\

\noindent{\bf Problem 10}. Does every polygon have a periodic orbit?\\

Every rational polygon has periodic orbits, and much is known
about them. Certain classes of irrational polygons have periodic
orbits \cite{Kolan,GuTr}. Every acute triangle has a classical
periodic orbit - the Fagnano orbit \cite{Gut4}. It corresponds to
the inscribed triangle of minimal perimeter. It is not known if
every acute triangle has other periodic orbits; it is also not
known if every obtuse triangle has a periodic orbit
\cite{GaStVor,Hunger}. See \S~~\ref{update} for updates and
elaborations.

A periodic orbit with an even number of segments is contained in a
{\em parallel band} of periodic orbits of the same length. The
boundary components of a band are concatenations of singular
orbits, the so-called {\em generalized diagonals} \cite{Kat}.
These are the billiard orbits with endpoints at the corners.
Periodic orbits with an odd number of segments (e. g., the Fagnano
orbit) are isolated. They seem to be rare; a rational polygon has
at most a finite number of them. Denote by $f_P(\ell)$ the number
of periodic bands of length at most $\ell$. This counting function
for periodic billiard orbits in polygons grows subexponentially
\cite{Kat,GuHa95,GuHa}. Conjecturally, there should be a universal
polynomial upper bound on $f_P(\cdot)$. See section~~\ref{update}
for elaborations.\\

\noindent {\bf Problem 11}. Find efficient upper and lower bounds
on $f_P$ for irrational polygons.\\

From now until the end of this section we consider only  rational
polygons. By results of Masur \cite{Mas1,Mas2} and Boshernitzan
\cite{Bosh1,Bosh2}, there exist numbers $0 < c_*(P) \le c^*(P)
<\infty$ such that $c_*(P)\ell^2 \le f_P(\ell) \le c^*(P)\ell^2$
for sufficiently large $\ell$. We will refer to these inequalities
as the {\em quadratic bounds on periodic billiard orbits}. The
numbers $c_*(P),c^*(P)$ are the {\em quadratic constants}.\\

\noindent {\bf Problem 12}. Find efficient estimates for quadratic constants.\\

In all known examples $f_P(\ell)/\ell^2$ has a limit, i. e.,
$c_*(P)=c^*(P)=c(P)$. In this case, we say that the {\em polygon
has quadratic asymptotics}. The preceding definitions and
questions have obvious counterparts for translation surfaces,
where periodic billiard orbits are replaced by closed geodesics.
There is a special class of polygons and surfaces: Those
satisfying the {\em lattice property}, or, simply, the {\em
lattice polygons and translation surfaces} \cite{Veech1,Guj,GuJ}.
They have quadratic asymptotics, and there are general expressions
for their quadratic constants \cite{Gut1,Veech1,Veech3,Guj,GuJ}.
Besides, the quadratic constants for several lattice polygons are
explicitly known \cite{Veech1,Veech3}. These formulas contain
rather subtle arithmetic. The (normalized) quadratic constant of
any lattice polygon is an algebraic number \cite{GuJ}.

There is a lot of information about lattice polygons and lattice
translation surfaces. Lattice polygons seem to be very rare. All
acute lattice triangles have been determined \cite{KenSmi,Pu}. For
obtuse triangles the question is still open. There is a
$\mac$-invariant Lebesgue-class finite measure on the {\em strata}
of the {\em moduli space} $\mm_g$ of translation surface of fixed
genus \cite{EskMas}. The strata correspond to the partitions
$2g-2=p_1+\cdots+p_t:p_i\in\N$. Hence, we can speak of a {\em
generic translation surface} in $\mm_g(p_1,\dots,p_t)$. The
quadratic constant of the generic translation surface
$S\in\mm_g(p_1,\dots,p_t)$ depends only on the stratum
\cite{EskMas,EskOk,EskMasZor}. These results crucially use the
relationship between the {\em Teichm\"uller flow} \cite{ImaTani}
on $\mm_g$ and linear flows on translation surfaces. See the
surveys \cite{MasTab,DeM,Avi} for this material.

The results for generic translation surfaces have no consequences
for rational polygons, since the generic translation surface does
not correspond to a polygon. However, an extension of the
Teichm\"uller flow approach establishes the quadratic asymptotics
for a special, but nontrivial class of rational polygons \cite{EskMasSchm}.\\

\noindent {\bf Problem 13}. Does every rational polygon have quadratic asymptotics?\\

The following section contains updates, elaborations, and
extensions of the preceding material.

\section{Comments, updates, and  extensions}\label{update}
At the time of writing this text, the problems discussed in the
survey \cite{Gut5} remain open. However, the works that have since
appeared contain substantial relevant material. The main purpose
of this section is to comment on this material. At the same time,
we take the opportunity to add a few extensions and ramifications
that for some reasons did not appear in the survey \cite{Gut5}.
Accordingly, we have updated and expanded the bibliography.

A few recent books contain discussions of the billiard ball
problem. The book \cite{Chernov4} is a thorough exposition of the
hyperbolic billiard dynamics. See \S~\ref{hyperbolic}. The
treatise \cite{Berger1} contains several discussions of
connections between the billiard ball problem and geometry. The
material in \cite{Berger1} is relevant for all three types of
billiard tables discussed above. The latter can be also said about
\cite{Tabach1}. However, the two books are fundamentally
different. The book \cite{Berger1} is a comprehensive treatise on
geometry, where the billiard ball problem is but one of a
multitude of illustrative examples and applications, while
\cite{Tabach1} is addressed primarily to young American students,
and discusses a few instances pertaining to the geometric aspect
of the billiard. The book \cite{Gruber} exposes several
applications of convex geometry to the billiard ball problem. See
the material in \S~\ref{elliptic}.

\medskip

\subsection{The Birkhoff conjecture and related material} \label{birkhoff_sub}
\hfill \break The traditional formulation of the Birkhoff
integrability conjecture is in terms of the billiard map. See
Definition~~\ref{integrable-def} and Problem $1$. The billiard
flow on a smooth, convex table $Y$ is a Hamiltonian system with
$2$ degrees of freedom. The Hamiltonian version of Problem $1$ is
as follows:

\medskip

\noindent {\bf Problem 1a}. Let $Y$ be a smooth, convex billiard
table. If the billiard flow on $Y$ is an integrable Hamiltonian system, then $Y$ is an ellipse.\\

Although Problem 1 and Problem 1a are obviously related, a
resolution of either one of them would not directly imply a
resolution of the other. The results of S. Bolotin \cite{Bolotin1}
provide evidence supporting a positive resolution of Problem 1a.

Let $\gamma\subset\RR$ be a closed convex curve. Denote by
$|\gamma|$ its perimeter. The {\em string construction} associates
with any $\ell>|\gamma|$ a  closed convex curve $G(\gamma,\ell)$
containing $\gamma$ in its interior. The curve $G(\gamma,\ell)$ is
obtained by the following ``physical process''. We take a ring of
length $\ell$ made from a soft, non-stretchable material and wrap
it around $\gamma$. We pull the ring tight with a pencil; then,
holding it tight, we rotate the pencil all the way around
$\gamma$. The moving pencil will then trace the curve
$G(\gamma,\ell)$. A gardener could use this process to design
fences around his flower beds. For this reason, the procedure is
sometimes called the {\em gardener construction} \cite{Berger1}.

Let $Y=Y(\gamma,\ell)\subset\RR$ be the billiard table whose
boundary is $G(\gamma,\ell)$. Then $\gamma$ is a caustic for the
billiard on $Y$ \cite{Laz,GutKat,Tabach1}. If $\gamma_1$ is an
ellipse and $\ell_1>|\gamma_1|$, then
$\gamma_2=G(\gamma_1,\ell_1)$ is a confocal ellipse. Let now
$\ell_2>|\gamma_2|$. Then $\gamma_3=G(\gamma_2,\ell_2)$ is a
confocal ellipse containing $\gamma_1$. There is a unique $\ell_3$
such that $\gamma_3=G(\gamma_1,\ell_3)$. This transitivity
property of gardener's construction is a consequence of the
integrability of billiard on an ellipse. The Birkhoff conjecture
suggests the following problem.

\medskip

\noindent {\bf Problem 1b}. Ellipses are the only closed convex
curves satisfying the above transitivity.\\

This entirely geometric variant of the Birkhoff conjecture is due
to R. Melrose. A positive solution of Problem 1b is the subject of
the PhD thesis of Melrose's student E. Amiran \cite{Amiran1}.
However, the work \cite{Amiran1} contains a serious gap, and the
question remains open.

Let now $\gamma$ be any closed convex curve. For $\ell>|\gamma|$
let $Y(\ell)=Y(\gamma,\ell)$ be the corresponding family of convex
billiard tables, and let $0<\rho_{\gamma}(\ell)<1/2$ be the
rotation number of the caustic $\gamma\subset Y(\ell)$. The {\em
rotation  function} $\rho(\ell)=\rho_{\gamma}(\ell)$ is continuous
and monotonically increasing, but not strictly, in general. Let
$r\in(0,1/2)$ be such that $\rho^{-1}(r)=[a(r),b(r)]$ is a
nontrivial interval. Then i) $r$ is rational; ii) for
$\ell\in[a(r),b(r)]$ the billiard map of $Y(\ell)$ restricted to
the caustic $\gamma$ is not a rotation. The converse also holds;
$[a(r),b(r)]$ are the {\em phase locking intervals}. The above
situation is a special case of the dynamical phenomenon called
{\em phase locking}. It is characteristic for one-parameter
deformations in elliptic dynamics. See \cite{GutKni} for a study
and a detailed  discussion of this phase locking when $\gamma$ is
a triangle.

Let now $\gamma$ be a convex polygon. Then the $C^1$ curve
$G(\gamma,\ell)$ is a concatenation of arcs of ellipses with foci
at the corners of $\gamma$. At the points of transition between
these elliptic arcs, typically, only one of the two foci changes,
causing a jump in the curvature. Thus, a typical billiard table,
say $Y$, obtained by this construction, is strictly convex,
piecewise analytic, but not $C^2$. The boundary $\bo Y$ contains a
finite number of points where the curvature jumps. A. Hubacher
studied billiard tables of this class \cite{Hubacher}. She proved
that there is an open neighborhood $\Omega\subset Y$ of $\bo Y$
such that any caustic in $Y$ belongs to the complement of $\Omega$
\cite{Hubacher}.

Recall that $\gamma$ is a caustic of $Y(\gamma,\ell)$ for any
$\ell>|\gamma|$. There is an analogy between Hubacher's theorem
and a result in \cite{GutKat} which says that billiard caustics
stay away from the table's boundary if it contains points of very
small curvature. This result is a quantitative version of Mather's
theorem \cite{Mather1} that insures nonexistence of caustics if
$\bo Y$ has points of zero curvature. Hubacher's theorem replaces
them with jump points of the curvature. As opposed to
\cite{GutKat}, the work \cite{Hubacher} does not estimate the size
of the {\em region $\Omega\subset Y$ free of caustics}. It is
plausible that for the typical table $Y(\ga,\ell)$ the free of
caustics region $\Omega=\Omega(\ga,\ell)$ is the annulus between
$\bo Y(\gamma,\ell)$ and $\gamma$.

There are polygons $\ga$ such that for special values of the
string length $\ell$ the boundary $\bo Y(\gamma,\ell)$ is a $C^2$
curve. Let $\ga$ be the regular hexagon with the unit side length.
Then $\bo Y(\gamma,14)$ is a $C^2$ curve. The work of H. Fetter
\cite{Fet2012} studies the billiard on $Y(\gamma,14)$. Fetter
suggests that the billiard on $Y(\gamma,14)$ is integrable, and
thus $Y(\gamma,14)$ is a counterexample to the Birkhoff
conjecture. However, the evidence of integrability of
$Y(\gamma,14)$ presented in \cite{Fet2012} is mostly numerical.
The present author believes that the further investigation of the
billiard on $Y(\gamma,14)$ will confirm the Birkhoff conjecture.

The subject of \cite{KalSor} is a billiard version of the famous
question of Marc Kac: {\em ``Can one hear the shape of a drum?}''
\cite{Kac}. As an application of their results, the authors prove
a {\em conditional version} of the Birkhoff conjecture.

\subsection{The Ivrii conjecture and related material} \label{ivrii_sub}
\hfill \break Investigations of Problem 3i) (i.e., the Ivrii
conjecture) and Problem 3ii) (i.e., the claims $I_n$ for $n>2$),
as well as related questions, remain active. The work
\cite{Baryshnikov1} develops a functional theoretic approach to
study billiard caustics. As a byproduct, \cite{Baryshnikov1}
contains yet another proof of Rychlik's theorem stating that the
set of $3$-periodic billiard orbits has measure zero. The preprint
\cite{Gut9} announced a solution of the Ivrii conjecture.
Unfortunately, the work \cite{Gut9} contains a mistake; thus the
conjecture remains open. The paper \cite{GluKudr} announced a
proof of the claim $I_4$: The set of $4$-periodic billiard orbits
has measure zero. The work \cite{GluKudr} states several
propositions implying the claim, and explains the strategy of
their proofs. The approach of \cite{GluKudr} is based on a study
of certain foliations, on one hand, and a very detailed analysis
of singularities of certain mappings, on the other hand. Complete
proofs should appear shortly. There is a curious connection
between the Ivrii conjecture and the subject of {\em invisibility}
\cite{PlakRosh}.

It is natural to investigate the counterparts of the Ivrii
conjecture for the billiard on (simply connected) surfaces of
constant curvature. The billiard on $\RR$ corresponds to the zero
curvature, $\kappa=0$. Multiplying a constant curvature by a
positive factor does not qualitatively change the geometry; thus,
it suffices to consider the two cases $\kappa=\pm 1$. The surfaces
in question are the hyperbolic plane, $\kappa=-1$, and the round
unit sphere, $\kappa=1$. Let us denote them by $\HH^2$ and
$\SSS^2$ respectively. On $\SSS^2$ the immediate analog of the
Ivrii conjecture fails. The paper \cite{GT} contains an example of
a (not strictly) smooth, convex billiard table in $\SSS^2$ with an
open set of periodic orbits. This observation shows the subtlety
of the Ivrii conjecture for $\RR$. The work \cite{BlKiNaZh}
contains a detailed study of $3$-periodic orbits for billiard
tables in $\HH^2$ and $\SSS^2$. It shows, in particular, that the
set of $3$-periodic billiard orbits on a Birkhoff billiard table
in $\HH^2$ has measure zero. This is the counterpart of Rychlik's
theorem for the hyperbolic plane.


\subsection{Extensions of the material in \S~\ref{elliptic}} \label{elliptic_sub}
\hfill \break Nontrivial billiard properties can be roughly
divided into three categories: i) Those that hold for all tables;
ii) Those that hold for the typical table; iii) Those that hold
for special billiard tables. Studies in category iii) can be
described as follows: Let P be a property satisfied by a very
particular billiard table, e. g., the round disc. Are there
non-round tables that have property P? If the answer is ``yes'',
then describe the billiard tables having property P.

The following example illustrates the situation. Let
$0<\alpha\le\pi/2$ be an angle. Let $Y\subset\RR$ be a Birkhoff
billiard table. We say that the table $Y$ has property
$P_{\alpha}$ if every chord in $Y$ that makes angle $\alpha$ with
$\partial Y$ at one end also makes angle $\alpha$ with $\partial
Y$ at the other end. The round table has property $P_{\alpha}$ for
any $\alpha$. A billiard table with property $P_{\alpha}$ has a
very special caustic $\Gamma_{\alpha}$; we will say that
$\Gamma_{\alpha}$ is a {\em constant angle caustic}. Let
$\rho(\het),0\le\het\le 2\pi,$ be the radius of curvature for
$\partial Y$. Tables with the caustic $\Gamma_{\pi/2}$ are well
known to geometers: Their boundaries are the {\em curves of
constant width} \cite{Berger1,Gruber}. A curve $\partial Y$ has
constant width if an only if its radius of curvature satisfies the
identity
$$
\rho(\het)+\rho(\het+\pi)=\const.
$$
In particular, there are non-round infinitely smooth, and even
real analytic billiard tables in this class.

For $0<\alpha<\pi/2$ let $\pp_{\alpha}$ be the class of non-round
tables with the property $P_{\alpha}$. The author has investigated
the class $\pp_{\alpha}$ about 20 years ago and reported the
results at the 1993 Pennsylvania State University Workshop on
Dynamics \cite{Gut+}. The main results in \cite{Gut+} are as
follows: The class $\pp_{\alpha}$ is nonempty if and only if
$\alpha$ satisfies
\begin{equation}  \label{alfa-eq}  
\tan(n\alpha)=n\tan\alpha
\end{equation}
for some $n>1$. The set $A_n\subset(0,\pi/2)$ of solutions of
equation~~\eqref{alfa-eq} is finite and nonempty for $n\ge 4$. For
every $\alpha\in A_n$ there is an analytic family of (nonround)
distinct, convex, real analytic tables
$Y_{\alpha}(s)\in\pp_{\alpha}:0<s<1$. As $s\to 0$, the tables
$Y_{\alpha}(s)$ converge to the unit disc. The limit
$Y_{\alpha}(1)$ of $Y_{\alpha}(s)$, as $s\to 1$, also exists, but
has corners.

It turns out that planar regions satisfying property $P_{\alpha}$
for some angle $\alpha$ are of interest in the mathematical fluid
mechanics. Besides the concept of Archimedean
floating,\footnote{In fact, this concept goes back to Aristotle.}
there is a concept of capillary {\em floating in neutral
equilibrium at a particular contact angle}. This concept goes back
to Thomas Young \cite{Young} and was further developed and
investigated by R. Finn \cite{Finn1}. If $Y\in\pp_{\alpha}$, then
the infinite cylinder $C=Y\times\R$ floats in neutral equilibrium
at the contact angle $\pi-\alpha$ {\em at every orientation}. The
work \cite{Gut7} is a detailed exposition of the results in
\cite{Gut+} aimed, in particular, at the mathematical fluid
mechanics readership. See \cite{Tabach3,Cyr} for related
investigations and \cite{Gut12} for additional comments.

\medskip

\subsection{Comments and updates for the material in \S~\ref{parabolic}} \label{parabolic_sub}
\hfill \break Although Problem 10 remains open, recent
publications \cite{Schw06,Schw09,HooSchw} provided substantial
evidence toward the positive answer, i.e., that every polygon does
have a periodic orbit. These papers investigate periodic billiard
orbits in obtuse triangles. Let $\Delta(\alpha,\beta,\ga)$ be the
triangle with the angles $\alpha\le\beta\le\ga$. In
\cite{Schw06,Schw09} R. Schwartz proves that if $\gamma$ is less
than or equal $100$ degrees, then $\Delta(\alpha,\beta,\ga)$ has a
periodic orbit. In \cite{HooSchw} Hooper and Schwartz prove that
if the angles $\alpha$ and $\beta$ are sufficiently close, then
$\Delta(\alpha,\beta,\ga)$ has a periodic orbit. If
$\alpha=\beta$, then $\Delta$ is an isosceles triangle. It is well
known and elementary that isosceles triangles have periodic
orbits. The main theorem in \cite{HooSchw} says that any triangle
which is sufficiently close to an isosceles one, has periodic
orbit. Note that in the Hooper-Schwartz theorem the angles
$\alpha$ and $\beta$ can be arbitrarily small, thus $\gamma$ can
be arbitrarily close to 180 degrees.

Let $\Delta$ be a triangle, and let $a,b,c$ be the sides of
$\Delta$. The works \cite{Schw06,Schw09,HooSchw} build on the
approach of \cite{GaStVor}, where periodic billiard orbits were
coded by words on the alphabet $\{a,b,c\}$. The paper
\cite{GaStVor} studied relationships between periodic orbits and
the associated words. A periodic billiard orbit is {\em stable} if
it persists under all sufficiently small deformations of $\Delta$.
By \cite{GaStVor}, the stability of an orbit is equivalent to a
combinatorial property of the associated word. For the sake of
brevity, I will simply say that the words on the alphabet
$\{a,b,c\}$ satisfying this property are {\em stable}. Let
$\ttt\subset\RR$ be the moduli space of triangles. Slightly
simplifying the situation, we assume that $\ttt=\{(x,y)\in\RR:0\le
x,y\le 1\}$. The subsets of obtuse (resp. isosceles) triangles are
given by $x+y<1$ (resp. $x=y$). Let $W_k$ be the set of stable
words of length $k$, and let $W=\cup_kW_k$. For $w\in W$ let
$\ttt_w\subset\ttt$ be the open set of triangles having a periodic
orbit with the code $w$. These are the {\em tiles} in the
terminology of \cite{HooSchw}.

The approach of Hooper and Schwartz is to exhibit a {\em
sufficient set} $W_{suff}\subset W$ so that the tiles
$\{\ttt_w:w\in W_{suff}\}$ cover the targeted part of $\ttt$. It
goes without saying that this idea cannot be implemented without
substantial computer power. The computer program ``MacBilliards''
created by Hooper and Schwartz does the job. Besides providing us
with ample evidence towards the conjecture that every polygon has
periodic billiard orbits, the works \cite{Schw06,Schw09,HooSchw}
have established several facts that show just how intricate  the
matter is. In some cases, every finite set is insufficient; then
\cite{HooSchw} sufficient infinite sets.

\medskip

Unfortunately, the survey \cite{Gut5} has omitted the subject of
{\em complexity of billiard orbits} in polygons, which is very
close to Problem 11. We will briefly discuss it below. Let $P$ be
a polygon, and let $\A=\{a,b,c,\ldots\}$ be the set of its sides.
Following a finite billiard orbit $\gamma$ and recording the sides
that it encounters, we obtain a word, $w(\gamma)$, on the alphabet
$\A$. We  say that $w(\gamma)$ is the {\em code} of $\gamma$; the
number of letters in $w(\gamma)$ is the {\em combinatorial length}
of $\gamma$. Let $W_n(P)$ be the set of codes of all billiard
orbits with combinatorial length $n$. The function $F(n)=|W_n(P)|$
is the {\em full complexity} of the billiard on $P$. Imposing
various restrictions on the billiard orbits, we obtain conditional
or {\em partial complexities} $F_*(n)$. For instance, the function
$f_P(n)$ in Problem 11 is the {\em periodic complexity} for the
billiard on $P$. Thus, Problem 11 is a special
case of the following open question.\\


\noindent {\bf Problem 14}. Find nontrivial, explicit bounds on
partial billiard complexities for the {\em general
polygon}.\footnote{The
designation {\em general polygon} may be replaced by {\em irrational polygon}.}\\

By \cite{Kat,GuHa95,GuHa}, the full complexity of the billiard in
any polygon is subexponential. This means that for $n$
sufficiently large $F(n)<e^{an}$ for any $a>0$. Note that this
result provides no subexponential bound on $F(n)$. Very few
nontrivial bounds on billiard complexities are known. Let $f_P(n)$
be the number of codes for periodic orbits in $P$ of length at
most $n$. For any $k\in\N$ Hooper constructed open sets $\x_k$ in
the moduli space of polygons such that for $P\in\x_k$ the function
$f_P(n)$ grows faster than $n\log^kn$ \cite{Hoop}. Let $g_P(n)$ be
the number of generalized diagonals in $P$ of length at most $n$.
Let $\ttt_3$ be the moduli space of triangles, endowed with the
Lebesgue measure. Scheglov \cite{Sche2012} showed that for almost
every $P\in\ttt_3$ and any $\ep>0$ the inequality
$g_p(n)<\const\exp(n^{\sqrt{3}-1+\ep})$ holds. Note that
Scheglov's result yields explicit subexponential upper bounds on
the full and the periodic complexities for almost every triangle.

To make Problem 14 more concrete, we state below a widely accepted
conjecture.

\medskip

\noindent {\bf Conjecture 1}. There is $d\ge 3$ such that the full
billiard  complexity for any polygon has a cubic lower bound and a
degree $d$ upper bound.\\

The claim is established only for rational polygons, with $d=3$.
This is a consequence of the results of Masur \cite{Mas1,Mas2}
about saddle connections on translation surfaces. We will now
briefly discuss recent results on partial complexities that
provide support for it \cite{GutRam}. Let $P$ be an arbitrary
polygon. Let $0\le\theta< 2\pi$ be a direction, and let $z\in P$
be a point. Coding those billiard orbits that start off in
direction $\theta$ (resp. from the point $z$) and counting the
number of words of length less than or equal to $n$, we obtain the
{\em directional complexity} $F_{\theta}(n)$ (resp. {\em position
complexity} $F_z(n)$) for the billiard on $P$. It was known that
the directional complexity grows polynomially, i.e., there is
$d>0$ depending only on $P$ such that $F_{\theta}(n)=O(n^d)$ for
all directions $\theta$  \cite{GuTr}. By \cite{GutRam}, for any
$P$ and any $\ep>0$ for almost all directions $\theta$, we have
$F_{\theta}(n)=O(n^{1+\ep})$. Another result in \cite{GutRam} says
that for any $P$ and any $\ep>0$, for almost all points $z\in P$
we have $F_z(n)=O(n^{2+\ep})$.\footnote{This establishes, in
particular, the unpublished results of Boshernitzan \cite{Bosh2}.}
The proofs are based on the relationships between the average
complexity and individual complexities. The concept of a {\em
piecewise convex polygon exchange} introduced in \cite{GT1} yields
a new approach to the billiard complexities. This approach works
for the polygonal billiard on surfaces of arbitrary constant
curvature.


\subsection{Ramifications and extensions of the polygonal billiard} \label{noncomp_sub}
\hfill\break In this section we briefly report on several
generalizations of the subject in section \S~\ref{parabolic} that
did not get discussed in \cite{Gut5}. The main direction in these
developments is to replace the usual polygons by noncompact or
infinite polygons. The work in this direction started already at
the end of the last century \cite{DegDelLen98,DegDelLen2000}, and
flourished after the turn of the century. The noncompact polygons
in \cite{DegDelLen98,DegDelLen2000} are semi-infinite stairways. A
stairway in these works is a noncompact polygon, say $P$, with
infinitely many vertical and horizontal sides, and of finite area.
In particular, $P$ is a {\em rational noncompact polygon}, and we
have the obvious family of directional billiard flows
$b_{\het}^t,0\le\het\le\pi/2,$ on $P$. In addition to studying the
ergodicity of these flows, the papers
\cite{DegDelLen98,DegDelLen2000} investigate {\em escaping orbits}
in $P$, a new phenomenon caused by the noncompactness of $P$.

There are many kinds of noncompact polygons, both rational and
irrational, bounded or unbounded, with finite or infinite area.
More generally, there are {\em noncompact polygonal surfaces}
\cite{Gut09}, and in particular, {\em noncompact translation
surfaces} \cite{Hoop08,Gut09}. In this {\em brave new world},
there are infinite (as opposed to semi-infinite) stairways
\cite{HooHubWei}, infinite coverings of compact polygonal surfaces
\cite{Gut09,FrUl11,CoGu}, and even fractals, e.g., the Koch
snowflake \cite{LapNie}. Besides being of interest on its own, the
billiard on noncompact polygons comes up in physics and
engineering \cite{GuKac,Gut85,Gut87,DanGu,BaKhMaPl}.

All of the questions concerning the dynamics and geometry for
compact polygons, translation surfaces, etc, have obvious
counterparts in the noncompact world.\footnote{Except for the
billiard on a fractal table where these counterparts are not
obvious \cite{LapNie}.} The answers to these questions are
sometimes unexpected \cite{FrUl11,CoGu,Trevino}, not at all
analogous to the answers in the compact world. Besides, there are
problems in the noncompact world that do not arise in the
classical setting. One of them is the conservativity of billiard
dynamics. In the classical setting the conservativity is ensured
by the Poincar\'e recurrence theorem. We conclude this necessarily
incomplete survey of noncompact polygonal billiard with a brief
discussion of the billiard for a classical family of noncompact,
doubly periodic polygons.

\medskip

Let $0<a, b <1$, and let $R(a,b)$ be the $a\times b$ rectangle.
Denote by $R_{(0,0)}$ the upright rectangle $R(a,b)$ centered at
$(1/2,1/2)$. For $(m,n)\in\ZZ$ set $R_{(m,n)}=R_{(0,0)}+(m,n)$.
The region
$\tP(a,b)=\RR\setminus\left(\cup_{(m,n)\in\ZZ}R_{(m,n)}\right)$ is
a noncompact rational polygon of infinite area. We will refer to
it as the {\em rectangular Lorenz gas}. The randomized version of
$\tP(a,b)$ is the famous {\em wind-tree model} of statistical
physics \cite{Eh}. For $0\le\het<\pi/2$ denote by $b_{\het}^t$ the
directional billiard flows on $\tP(a,b)$. The work \cite{CoGu}
works out the ergodic decomposition of $b_{\het}^t$ on $\tP(a,b)$
for particular directions $\het$ and sufficiently small $a$ and
$b$, provided that $a/b$ be irrational. The directions in question
are $\het=\arctan(q/p)$ corresponding to $(p,q)\in\NN$ such that
the flow line of $b_{\het}^t$ emanating from a corner of
$R_{(0,0)}$ reaches the homologous corner of $R_{(p,q)}$ bypassing
the obstacles. Although these directions form a finite set, it is
asymptotically dense as $a,b$ go to zero. The dissipative
component of $b_{\het}^t$ is spanned by the straight
lines\footnote{With the slopes $\pm q/p$.} avoiding the obstacles.
The conservative part of $b_{\het}^t$ decomposes as a direct sum
of $2pq$ identical ergodic flows. The decomposition is as follows.
There is a certain subgroup $H_{(p,q)}\subset\ZZ$ of index $2pq$.
Let $S_1,\dots,S_{2pq}\subset\ZZ$ be its cosets. The ergodic
component of $b_{\het}^t$ corresponding to the coset $S_i,1\le
i\le 2pq,$ is spanned by the billiard orbits that encounter the
obstacles $R_{(m,n)}$, where $(m,n)\in S_i$. Thus, the ergodic
decomposition of the flow in the $(p,q)$ direction is induced by a
natural partition of the set of obstacles in the configuration
space. For instance, the two ergodic components of the flow in the
direction $\pi/4$ correspond to the billiard orbits on $\tP(a,b)$
that encounter the obstacles $R_{(m,n)}$ with even and odd  $m+n$
respectively.

These results follow from the ergodicity of certain $\ZZ$-valued
cocycles over irrational rotations \cite{CoGu} established by the
classical methods of ergodic theory for infinite invariant
measures \cite{Ar,Co,Schm}. The explicit ergodic decomposition of
the (conservative part of) the directional flows yields nontrivial
consequences for the recurrence and the spatial distribution of
typical billiard orbits \cite{CoGu}. The recurrence in the
direction $\pi/4$ has been discussed in the physics literature
\cite{HW80}. Judging by the results in \cite{CoGu,FrUl11} etc, the
noncompact polygonal billiard may yield further surprises.


\subsection{Security for billiard tables and related questions} \label{security_sub}
\hfill \break This subject arose quite recently. It has to do with
the geometry of billiard orbits as curves on the configuration
space. Let $\Om\subset\RR$ be any billiard table. For any pair
$x,y\in\Om$ (in particular for pairs $x,x$) let $\Ga(x,y)$ be the
set of billiard orbits in $\Om$ that connect $x$ and $y$. We say
that the {\em pair $x,y$ is secure} if there exists a finite set
$\{z_1,\ldots,z_n\}\subset\Om\setminus\{x,y\}$ such that every
$\ga\in\Ga(x,y)$ passes through some $z_i$. We say that
$z_1,\ldots,z_n$ are {\em blocking points} for $x,y$. We call
$\Om$ secure if every pair of points is secure. Thus, to show that
$\Om$ is {\em insecure} means to find a pair $x,y$ that cannot be
blocked with a finite set of blocking points. The questions that
arise are: i) What tables $\Om$ are secure?; ii) If $\Om$ is
insecure, how insecure is it? For instance, is it true that every
pair $x,y\in\Om$ is insecure, that almost all pairs $x,y\in\Om$
are insecure, etc.

The subject, in disguise of {\em problems about bodyguards} came
up in the Mathematical Olympiad literature. The recent interest in
security got triggered by \cite{HS}. The authors, who were then
students at Cambridge University, studied the security of
polygons. The main claim in \cite{HS} is that every rational
polygon is secure. Regular $n$-gons provide a counterexample to
the claim: A regular $n$-gon is secure if and only if $n=3,4,6$
\cite{Gut6}. See \cite{Gut7} for related results. Although the
statement is elementary, the proof is not. The claim follows from
a study of security in translation surfaces, and is based on
\cite{Gut1,Veech1,Veech3,Guj,GuJ,GuHuSc}. However, the general
study of security for polygons has just begun.

\medskip

\noindent {\bf Problem 15}. To characterize secure polygons. In
particular, establish a criterion of security for rational polygons.\\

Triangles with the angles of $30,60,90$ and $45,45,90$ degrees, as
well as the equilateral triangle and the rectangles are the only
polygons whose translation surfaces are flat tori, and hence their
billiard flows are integrable \cite{Gut1}. Following \cite{Gut1},
we will call them {\em integrable polygons}. A polygon $P$ is {\em
almost integrable} if it is tiled under reflections by one of the
integrable polygons \cite{GutKat1}. By \cite{Gut6}, every almost
integrable polygon is secure.\\


\noindent {\bf Conjecture 2}. A polygon is secure if and only if
it is almost integrable.\\

Almost integrable polygons are certainly rational. Conjecture 2
would imply, in particular, that every irrational polygon is
insecure. At present, the problem of insecurity for general
polygons seems hopeless. However, the security framework makes
sense for arbitrary billiard tables, and more generally, for
arbitrary riemannian manifolds (with boundary, corners, and
singularities, in general). Security of riemannian manifolds is
related to the growth of the number of connecting geodesics, and
hence to their {\em topological entropy} \cite{BuGut,Gut8}. Not
much is known about the security in non-polygonal billiard tables
$\Om$. Let $\Om$ be a Birkhoff billiard table. Approximating
$\bo\Om$ locally by the arcs of its osculating circles, it is
intuitively clear that pairs of sufficiently close points
$x,y\in\bo\Om$ are insecure; the work \cite{Tabach2} confirms it.

Much more is known about the security of compact, smooth
Riemannian manifolds. Flat manifolds are secure \cite{Gut6,BuGut}.
For a flat torus, the proof is  elementary \cite{Gut6}. For
general flat manifolds, this follows from the Bieberbach theorem
\cite{BuGut}.

\medskip

\noindent {\bf Conjecture 3}. A compact, smooth riemannian
manifolds is secure if and only if it is flat.\\

Conjecture 3 has been established for various classes of compact
Riemannian manifolds \cite{BuGut,GutSch,BaGu}. In particular, it
holds for all compact surfaces of genus greater than zero
\cite{BaGu}. Thus, among surfaces, it remains to prove it for
arbitrary smooth Riemannian metrics on the two-sphere. The work
\cite{GerbLi} gives examples of totally insecure real analytic
metrics on the two-sphere; \cite{GerbKu} shows that higher
dimensional compact Riemannian manifolds are generically insecure.

\medskip

\noindent{\bf Acknowledgements}. The author gratefully
acknowledges discussions with many mathematicians concerning this
survey and the comments of anonymous referees. The work was
partially supported by the MNiSzW Grant N N201 384834.

\end{document}

%% file: billi1.pstex_t
\begin{picture}(0,0)%
\includegraphics{billi1.pstex}%
\end{picture}%
\setlength{\unitlength}{1973sp}%
\begingroup\makeatletter\ifx\SetFigFont\undefined%
\gdef\SetFigFont#1#2#3#4#5{%
  \reset@font\fontsize{#1}{#2pt}%
  \fontfamily{#3}\fontseries{#4}\fontshape{#5}%
  \selectfont}%
\fi\endgroup%
\begin{picture}(8991,9140)(2049,-8461)
\put(5701,-8011){\makebox(0,0)[lb]{\smash{\SetFigFont{9}{10.8}{\rmdefault}{\mddefault}{\updefault}{\color[rgb]{0,0,0}$\theta$}%
}}}
\put(4501,-2911){\makebox(0,0)[lb]{\smash{\SetFigFont{12}{14.4}{\rmdefault}{\mddefault}{\updefault}{\color[rgb]{0,0,0}$Y$}%
}}}
\put(9226,-2536){\makebox(0,0)[lb]{\smash{\SetFigFont{12}{14.4}{\rmdefault}{\mddefault}{\updefault}{\color[rgb]{0,0,0}$\bo Y$}%
}}}
\put(4876,-8086){\makebox(0,0)[lb]{\smash{\SetFigFont{9}{10.8}{\rmdefault}{\mddefault}{\updefault}{\color[rgb]{0,0,0}$\theta$}%
}}}
\put(5251,-8461){\makebox(0,0)[lb]{\smash{\SetFigFont{10}{12.0}{\rmdefault}{\mddefault}{\updefault}{\color[rgb]{0,0,0}$x$}%
}}}
\put(4126,-6886){\makebox(0,0)[lb]{\smash{\SetFigFont{9}{10.8}{\rmdefault}{\mddefault}{\updefault}{\color[rgb]{0,0,0}$v$}%
}}}
\put(7126,-6136){\makebox(0,0)[lb]{\smash{\SetFigFont{9}{10.8}{\rmdefault}{\mddefault}{\updefault}{\color[rgb]{0,0,0}$v'$}%
}}}
\put(6601,-6661){\makebox(0,0)[lb]{\smash{\SetFigFont{9}{10.8}{\rmdefault}{\mddefault}{\updefault}{\color[rgb]{0,0,0}$y'$}%
}}}
\put(3301,-6886){\makebox(0,0)[lb]{\smash{\SetFigFont{9}{10.8}{\rmdefault}{\mddefault}{\updefault}{\color[rgb]{0,0,0}$y$}%
}}}
\end{picture}

%% file: billi2.pstex_t
\begin{picture}(0,0)%
\includegraphics{billi2.pstex}%
\end{picture}%
\setlength{\unitlength}{1973sp}%
\begingroup\makeatletter\ifx\SetFigFont\undefined%
\gdef\SetFigFont#1#2#3#4#5{%
  \reset@font\fontsize{#1}{#2pt}%
  \fontfamily{#3}\fontseries{#4}\fontshape{#5}%
  \selectfont}%
\fi\endgroup%
\begin{picture}(9973,7077)(1578,-7709)
\put(5026,-7636){\makebox(0,0)[lb]{\smash{\SetFigFont{9}{10.8}{\rmdefault}{\mddefault}{\updefault}{\color[rgb]{0,0,0}$x$}%
}}}
\put(11551,-4786){\makebox(0,0)[lb]{\smash{\SetFigFont{9}{10.8}{\rmdefault}{\mddefault}{\updefault}{\color[rgb]{0,0,0}$x_1$}%
}}}
\put(5551,-7411){\makebox(0,0)[lb]{\smash{\SetFigFont{9}{10.8}{\rmdefault}{\mddefault}{\updefault}{\color[rgb]{0,0,0}$\theta$}%
}}}
\put(11101,-5011){\makebox(0,0)[lb]{\smash{\SetFigFont{9}{10.8}{\rmdefault}{\mddefault}{\updefault}{\color[rgb]{0,0,0}$\theta_1$}%
}}}
\put(11251,-4486){\makebox(0,0)[lb]{\smash{\SetFigFont{9}{10.8}{\rmdefault}{\mddefault}{\updefault}{\color[rgb]{0,0,0}$\theta_1$}%
}}}
\put(4726,-2161){\makebox(0,0)[lb]{\smash{\SetFigFont{12}{14.4}{\rmdefault}{\mddefault}{\updefault}{\color[rgb]{0,0,0}$\bo Y$}%
}}}
\end{picture}

%% file: billi3.pstex_t
\begin{picture}(0,0)%
\includegraphics{billi3.pstex}%
\end{picture}%
\setlength{\unitlength}{1973sp}%
\begingroup\makeatletter\ifx\SetFigFont\undefined%
\gdef\SetFigFont#1#2#3#4#5{%
  \reset@font\fontsize{#1}{#2pt}%
  \fontfamily{#3}\fontseries{#4}\fontshape{#5}%
  \selectfont}%
\fi\endgroup%
\begin{picture}(10741,7599)(1268,-7798)
\end{picture}

%% file: billi4.pstex_t
\begin{picture}(0,0)%
\includegraphics{billi4.pstex}%
\end{picture}%
\setlength{\unitlength}{1973sp}%
\begingroup\makeatletter\ifx\SetFigFont\undefined%
\gdef\SetFigFont#1#2#3#4#5{%
  \reset@font\fontsize{#1}{#2pt}%
  \fontfamily{#3}\fontseries{#4}\fontshape{#5}%
  \selectfont}%
\fi\endgroup%
\begin{picture}(9558,9645)(1411,-8956)
\put(6136,389){\makebox(0,0)[lb]{\smash{\SetFigFont{12}{14.4}{\rmdefault}{\mddefault}{\updefault}{\color[rgb]{0,0,0}$D$}%
}}}
\put(1411,389){\makebox(0,0)[lb]{\smash{\SetFigFont{12}{14.4}{\rmdefault}{\mddefault}{\updefault}{\color[rgb]{0,0,0}$E$}%
}}}
\put(1456,-8956){\makebox(0,0)[lb]{\smash{\SetFigFont{12}{14.4}{\rmdefault}{\mddefault}{\updefault}{\color[rgb]{0,0,0}$A$}%
}}}
\put(10906,-4411){\makebox(0,0)[lb]{\smash{\SetFigFont{12}{14.4}{\rmdefault}{\mddefault}{\updefault}{\color[rgb]{0,0,0}$C$}%
}}}
\put(10891,-8911){\makebox(0,0)[lb]{\smash{\SetFigFont{12}{14.4}{\rmdefault}{\mddefault}{\updefault}{\color[rgb]{0,0,0}$B$}%
}}}
\end{picture}

%% file: billi5.pstex_t
\begin{picture}(0,0)%
\includegraphics{billi5.pstex}%
\end{picture}%
\setlength{\unitlength}{1973sp}%
\begingroup\makeatletter\ifx\SetFigFont\undefined%
\gdef\SetFigFont#1#2#3#4#5{%
  \reset@font\fontsize{#1}{#2pt}%
  \fontfamily{#3}\fontseries{#4}\fontshape{#5}%
  \selectfont}%
\fi\endgroup%
\begin{picture}(8609,8012)(1914,-8252)
\put(2446,-3946){\makebox(0,0)[lb]{\smash{\SetFigFont{10}{12.0}{\rmdefault}{\mddefault}{\updefault}{\color[rgb]{0,0,0}$m_1$}%
}}}
\put(3121,-6451){\makebox(0,0)[lb]{\smash{\SetFigFont{10}{12.0}{\rmdefault}{\mddefault}{\updefault}{\color[rgb]{0,0,0}$m_2$}%
}}}
\put(9481,-3751){\makebox(0,0)[lb]{\smash{\SetFigFont{10}{12.0}{\rmdefault}{\mddefault}{\updefault}{\color[rgb]{0,0,0}$m_3$}%
}}}
\end{picture}